\documentclass[10pt]{amsart}
\usepackage{palatino}
\usepackage[all]{xy,xypic}

\usepackage{amsfonts,amsmath,oldgerm,amssymb,amscd}

\newcommand{\ra}{\rightarrow}		

\newcommand{\by}[1]{\stackrel{#1}{\ra}}

\newcommand{\ol}{\overline}		
\newcommand{\wt}{\widetilde}
\newcommand{\iso}{\by \sim}

\newtheorem{theorem}{Theorem}[section]
\newtheorem{proposition}[theorem]{Proposition}
\newtheorem{lemma}[theorem]{Lemma}
\newtheorem{definition}[theorem]{Definition}
\newtheorem{corollary}[theorem]{Corollary}
\newtheorem{question}[theorem]{Question}

	\newcommand{\bz}{\mbox{$\mathbb Z$}}

\newcommand{\CC}{\mbox{$\mathcal C$}}

\newcommand{\CM}{\mbox{$\mathcal M$}}

	\newcommand{\p}{\mbox{$\mathfrak p$}}

\newcommand{\ot}{\mbox{\,$\otimes$\,}}	
\newcommand{\op}{\mbox{$\oplus$}}
\newcommand{\Spec}{\text{Spec}}

\newcommand{\Aut}{\mbox{\rm Aut\,}}	
\newcommand{\Hom}{\text{Hom}\,}

\newcommand{\Um}{\text{Um}}

\newcommand{\bp}{\begin{proposition}}
\newcommand{\ep}{\end{proposition}}
\newcommand{\bl}{\begin{lemma}}
\newcommand{\el}{\end{lemma}}
\newcommand{\bt}{\begin{theorem}}
\newcommand{\et}{\end{theorem}}
\newcommand{\bc}{\begin{corollary}}
\newcommand{\ec}{\end{corollary}}
\newcommand{\bd}{\begin{definition}}
\newcommand{\ed}{\end{definition}}

\def\rmk{\refstepcounter{theorem}\paragraph{{\bf Remark} \thetheorem}}

\newcommand{\wi}{\widetilde}

\def\rmk{\refstepcounter{theorem}\paragraph{{\bf Remark} \thetheorem}}

\def\proof{\paragraph{Proof}}

\def\definition{\refstepcounter{theorem}\paragraph{{\bf Definition} \thetheorem}}
\newcommand{\remark}{\rmk}

\title{Unimodular Elements in Projective
  Modules and an analogue of a result of Mandal} \author{Manoj
  K. Keshari} \address{Department of Mathematics, Indian Institute of
  Technology Bombay, Powai, Mumbai 400076, India.}
\email{keshari@math.iitb.ac.in} \author{Md. Ali Zinna}
\address{Department of Mathematics, Indian Institute of Technology
  Bombay, Powai, Mumbai, India - 400 076}
\email{zinna@math.iitb.ac.in} \subjclass[2000]{13C10}

\begin{document}
\maketitle
\section{Introduction}

{\it Throughout the paper, rings are commutative Noetherian and projective
modules are finitely generated and of constant rank.}

If $R$ is a ring of dimension $n$, then Serre \cite{Se} proved that
projective $R$-modules of rank $>n$ contain a unimodular element.
Plumstead \cite{p} generalized this result and proved that projective
$R[X]=R[\bz_+]$-modules of rank $>n$ contain a unimodular element.  Bhatwadekar
and Roy \cite{BRoy} generalized this result and proved that projective
$R[X_1,\ldots,X_r]=R[\bz_+^r]$-modules of rank $>n$ contain a unimodular element.

In another direction, if $A$ is a ring such that $R[X]\subset A\subset
R[X,X^{-1}]$, then Bhatwadekar and Roy \cite{br} proved that
projective $A$-modules of rank $>n$ contain a unimodular element.  Rao
\cite{Ra} improved this result and proved that if $B$ is a birational
overring of $R[X]$, i.e. $R[X]\subset B\subset S^{-1}R[X]$, where $S$
is the set of non-zerodivisors of $R[X]$, then projective $B$-modules
of rank $>n$ contain a unimodular element.  Bhatwadekar, Lindel and
Rao \cite[Theorem 5.1, Remark 5.3]{blr} generalized this result and
proved that projective $B[\bz_+^r]$-modules of rank $>n$ contain a
unimodular element when $B$ is seminormal. Bhatwadekar \cite[Theorem
  3.5]{Bh} removed the hypothesis of seminormality used in \cite{blr}.

All the above results are best possible in the sense that projective
modules of rank $n$ over above rings need not have a unimodular
element. So it is natural to look for obstructions for a projective
module of rank $n$ over above rings to contain a unimodular element.
We will prove some results in this direction.

Let $P$ be a projective $R[\bz_+^r][T]$-module of rank $n=\dim R$ such
that $P_f$ and $P/TP$ contain unimodular elements for some monic
polynomial $f$ in the variable $T$. Then $P$ contains a unimodular
element. The proof of this result is implicit in \cite[Theorem
  5.1]{blr}.  We will generalize this result to projective
$R[M][T]$-modules of rank $n$, where $M\subset \bz_+^r$ is a
$\Phi$-simplicial monoid in the class $\CC(\Phi)$. For this we need
the following result whose proof is similar to \cite[Theorem
  5.1]{blr}.

\begin{proposition}\label{roit1}
Let $R$ be a  ring and $P$
be a projective $R[X]$-module. Let $J\subset R$ be an ideal such that
$P_s$ is extended from $R_s$ for every $s\in J$. Suppose that

(a) $P/JP$ contains a unimodular element.

(b) If $I$ is an ideal of $(R/J)[X]$ of height $\text{rank}(P)-1$,
then there exist $\ol \sigma\in \Aut((R/J)[X])$ with $\ol \sigma(X)=X$
and $\sigma\in \Aut(R[X])$ with $\sigma(X)=X$ which is a lift of $\ol
\sigma$ such that $\ol \sigma(I)$ contains a monic polynomial in
the variable $X$.

(c) $EL(P/(X,J)P)$ acts transitively on $\Um(P/(X,J)P)$.

(d) There exists a monic polynomial $f\in R[X]$ such that $P_f$ 
contains a unimodular element.

Then the natural map $\Um(P)\ra \Um(P/XP)$ is surjective.
In particular, if $P/XP$ contains a unimodular element,
then $P$ contains a unimodular element.  
\end{proposition}

We prove the following result as an application of (\ref{roit1}).

\begin{theorem}\label{blr20}
Let $R$ be a ring of dimension $n$ and $M\subset \bz_+^r$ a
$\Phi$-simplicial monoid in the class $\CC(\Phi)$. Let $P$ be a
projective $R[M][T]$-module of rank $n$ whose determinant is extended
from $R$. Assume $P/TP$ and $P_f$ contain unimodular elements for
some monic polynomial $f$ in the variable $T$. Then the natural map $\Um(P)\ra
\Um(P/TP)$ is surjective. In particular, $P$ contains a unimodular
element.
\end{theorem}

Let $R$ be a ring containing $\mathbb{Q}$ of dimension $n\geq 2$. If
$P$ is a projective $R[X]$-module of rank $n$, then Das and Zinna
\cite{dz} have obtained an obstruction for $P$ to have a unimodular
element. Let us fix an isomorphism $\chi : L\iso \wedge^n P$, where $L$
is the determinant of $P$. To the pair $(P,\chi)$, they associated an element
$e(P,\chi)$ of the Euler class group $E(R[X],L)$ and 
proved that $P$ has a unimodular element if and only if $e(P,\chi)=0$
in $E(R[X],L)$ \cite{dz}. 

It is desirable to have such an obstruction for projective
$R[X,Y]$-module $P$ of rank $n$.  As an application of (\ref{blr20}),
we obtain such a result.  Recall that $R(X)$ denotes the ring obtained
from $R[X]$ by inverting all monic polynomials in $X$. Let $L$ be the
determinant of $P$ and $\chi:L\iso \wedge^n(P)$ be an isomorphism.  We
define the Euler class group $E(R[X,Y],L)$ of $R[X,Y]$ as the product
of Euler class groups $E(R(X)[Y],L\otimes R(X)[Y])$ of $R(X)[Y]$ and
$E(R[Y],L\otimes R[Y])$ of $R[Y]$ defined by Das and Zinna
\cite{dz}. To the pair $(P,\chi)$, we associate an element $e(P,\chi)$
in $E(R[X,Y],L)$ and prove the following result (\ref{obstruct}).

\begin{theorem}\label{obst}
Let the notations be as above. Then $e(P,\chi)=0$ in
$E(R[X,Y],L)$ if and only if $P$ has a unimodular element.
\end{theorem}

Let $R$ be a  local ring and $P$ be a projective
$R[T]$-module. Roitman \cite[Lemma 10]{ro} proved that
if the projective $R[T]_f$-module $P_f$ contains a unimodular element for some
monic polynomial $f\in R[T]$, then $P$ contains a unimodular element.
Roy \cite[Theorem
  1.1]{ry} generalized this result and proved that
if $P$ and $Q$ are projective $R[T]$-modules
with $\text{rank}(Q) < \text{rank}(P)$
such that $Q_f$ is a direct summand of $P_f$ for some
monic polynomial $f\in R[T]$, then $Q$ is a direct summand of $P$.
Mandal \cite[Theorem 2.1]{m} extended Roy's result to Laurent
polynomial rings.

We prove the following result (\ref{roy}) which gives Mandal's \cite{m}
in case $A=R[X,X^{-1}]$.  Recall that a monic polynomial $f\in R[X]$
is called \emph{special monic} if $f(0)=1$.

\bt\label{rm} 
Let $R$ be a  local ring and $R[X]\subset
A\subset R[X,X^{-1}]$. 
Let $P$ and $Q$ be two projective
$A$-modules with $\text{rank}(Q)< \text{rank}(P)$. If $Q_f$ is a
direct summand of $P_f$ for some special monic polynomial $f\in R[
X]$,
then $Q$ is also a direct summand of $P$.  \et

\section{Preliminaries}

\begin{definition} Let $R$ be a ring and $P$ be a projective $R$-module. An
element $p\in P$ is called \emph{unimodular} if there is a surjective
R-linear map $\varphi:P\twoheadrightarrow R$ such that
$\varphi(p)=1$. Note that $P$ has a unimodular element if and only if
$P\simeq Q \oplus R$ for some $R$-module $Q$.  The set of all
unimodular elements of $P$ is denoted by $\Um(P)$.
\end{definition}
\smallskip

\begin{definition}
Let $M$ be a finitely generated submonoid of $\bz_+^r$ of rank $r$
such that $M\subset \bz_+^r$ is an integral extension, i.e. for any
$x\in \bz_+^r$, $nx\in M$ for some integer $n>0$. Such a monoid $M$ is
called a {\it $\Phi$-simplicial} monoid of rank $r$ \cite{G2}.
\end{definition}
\smallskip

\begin{definition}
Let $M\subset \bz_+^{r}$ be a $\Phi$-simplicial monoid of rank $r$.
We say that $M$ belongs to the class $\CC(\Phi)$ if $M$ is seminormal
(i.e. if $x\in gp(M)$ and $x^2,x^3\in M$, then $x\in M$) and if we
write $\bz_+^r=\{t_1^{s_1}\ldots t_{r}^{s_r}|\, s_i\geq 0\}$, then for
$1\leq m \leq r$, $M_m=M\cap \{t_1^{s_1}\ldots t_{m}^{s_m}|\, s_i\geq
0\}$ satisfies the following properties: Given a positive integer $c$,
there exist integers $c_i > c$ for $i=1,\ldots, m-1$ such that for any
ring $R$, the automorphism $\eta \in Aut_{R[t_m]}(
R[t_1,\ldots,t_{m}])$ defined by $\eta(t_i)=t_i+t_m^{c_i}$ for $
i=1,\ldots, m-1$, restricts to an $R$-automorphism of $R[M_m]$.  It is
easy to see that $M_m\in \CC(\Phi)$ and rank $M_m=m$ for $1\leq m\leq
r$.
\end{definition}
\smallskip

\begin{example} The following monoids belong to 
$\CC(\Phi)$ \cite[Example 3.5, 3.9, 3.10]{KS}.

(i) If $M\subset \bz_+^2$ is a finitely generated and normal monoid
  (i.e. $x\in gp(M)$ and $x^n\in M$ for some $n>1$, then $x\in M$) of
  rank $2$, then $M\in \CC(\Phi)$.

(ii) For a fixed integer $n>0$, if
$M\subset \bz_+^r$ is the monoid generated by all monomials in
$t_1,\ldots,t_r$ of total degree $n$, then $M$ is a normal monoid of
rank $r$ and $M\in \CC(\Phi)$. In particular,
$\bz_+^r\in\CC(\Phi)$ and
$<t_1^2,t_2^2,t_3^2,t_1t_2,t_1t_3,t_2t_3> \in\CC(\Phi)$. 

(iii) The submonoid $M$ of $\bz_+^3$ generated by
$<t_1^2,t_2^2,t_3^2,t_1t_3,t_2t_3> \in\CC(\Phi)$.
\end{example}
\smallskip

\begin{remark}\label{rem1}
Let $R$ be a ring and $M\subset \bz_+^r=\{t_1^{m_1}\ldots
t_r^{m_r}|m_i \geq 0\}$ be a monoid of rank $r$ in the class
$\CC(\Phi)$. Let $I$ be an ideal of $R[M]$ of height $>\dim R$. Then by
\cite[Lemma 6.5]{G2} and \cite[Lemma 3.1]{KS}, there exists an
$R$-automorphism $\sigma$ of $R[M]$ such that $\sigma(t_r)=t_r$ and
$\sigma(I)$ contains a monic polynomial in $t_r$ with coefficients in
$R[M]\cap R[t_1,\ldots,t_{r-1}]$.
\end{remark}

We will state some results for later use.
\begin{theorem}\label{ks1}
\cite[Theorem 3.4]{KS}
Let $R$ be a  ring and $M$ be a $\Phi$-simplicial monoid
such that $M\in \CC(\Phi)$. Let $P$ be a projective $R[M]$-module of
rank $>\dim R$. Then $P$ has a unimodular element.
\end{theorem}

\begin{theorem}\label{dk1}
\cite[Theorem 4.5]{DK}
Let $R$ be a  ring and $M$ be a $\Phi$-simplicial
monoid. Let $P$ be a projective $R[M]$-module of rank $\geq max\{\dim
R+1,2\}$. Then $EL(P\oplus R[M])$ acts transitively on $\Um(P\oplus
R[M])$.
\end{theorem}

The following result is proved in \cite[Criterion-1 and Remark]{blr}
in case $J=Q(P,R_0)$ is the Quillen ideal of $P$ in $R_0$. The same proof works
in our case.

\begin{theorem}\label{criterion}
Let $R=\op_{i\geq 0}R_i$ be a graded ring and $P$ be a projective
$R$-module. Let $J$ be an ideal of $R_0$ such that $J$ is contained in
the Quillen ideal $Q(P,R_0)$. Let $p\in P$ be such that $p_{1+R^+}\in
\Um(P_{1+R^+})$ and $p_{1+J}\in \Um(P_{1+J})$, where $R^+=\op_{i\geq
  1}R_i$. Then $P$ contains a unimodular element $p_1$ such that
$p=p_1$ modulo $R^+P$.
\end{theorem}

The following result is a consequence of Eisenbud-Evans \cite{ee}, as stated in
\cite[p. 1420]{p}.

\begin{lemma}\label{EE}
Let $A$ be a ring and $P$ be a projective $A$-module of rank $n$. Let 
$(\alpha, a)\in (P^\ast \oplus A)$. Then there exists an element 
$\beta \in P^\ast$ such that {\rm ht}\,$(I_a)\geq n$,
where $I = (\alpha + a\beta)(P)$. In particular, if the ideal 
$(\alpha(P),a)$ has height $\geq n$, then {\rm ht}\,$I\geq n$. Further, if 
$(\alpha(P),a)$ is an ideal of height $\geq n$ and $I$ is a proper ideal of 
$A$, then {\rm ht}\,$I=n$.  
\end{lemma}

\section{Proofs of (\ref{roit1}), (\ref{blr20}) and (\ref{obst})}

\subsection{Proof of Proposition \ref{roit1}}

Let $p_0\in \Um(P/JP)$ and $p_1\in \Um(P/XP)$. Let $\wt {p_0}$ and
$\wt {p_1}$ be the images of $p_0$ and $p_1$ in $P/(X,J)P$.  By
hypothesis $(c)$, there exist $\wt \delta \in EL(P/(X,J)P)$ such that
$\wt \delta(\wt {p_0})=\wt {p_1}$. By \cite[Proposition 4.1]{BRoy}, $\wt \delta$
can be lifted to an automorphism $\delta$ of $P/JP$.  Consider the fiber product diagram for rings and modules
$$\xymatrix{ \frac {R[X]}{(XJ)} \ar@{->}[r] \ar@{->}[d] & \frac RJ[X] 
\ar@{->}[d] \\ \frac {R[X]}{(X)}
  \ar@{->}[r] &  \frac {R[X]}{(X,J)} \,.  } \hspace*{1in} 
\xymatrix{ \frac {P}{(XJ)P} \ar@{->}[r] \ar@{->}[d] & \frac P{JP} 
\ar@{->}[d] \\ \frac {P}{XP}
  \ar@{->}[r] &  \frac {P}{(X,J)P} \,.  } $$

We can patch
$\delta(p_0)$ and $p_1$ to get a unimodular element $p\in
\Um(P/XJP)$ such that $p=\delta(p_0)$ modulo $JP$ and $p=p_1$
module $XP$. Writing $\delta(p_0)$ by $p_0$, we assume that
$p=p_0$ modulo $JP$ and $p=p_1$ module $XP$.

Using hypothesis $(d)$, we get an element $q\in P$ such that the order
ideal $O_P(q)=\{\phi(q)|\phi\in \Hom_{R[X]}(P,R[X])\}$ contains a
power of $f$. We may assume that $f\in O_P(q)$.

Let ``bar'' denote reduction modulo the ideal $(J)$. Write $\ol P=\ol
{R[X]}p_0\oplus Q$ for some projective $\ol{R[X]}$-module $Q$ and
$\ol q=(\ol ap_0,q')$ for some $q'\in Q$. By Eisenbud-Evans (\ref{EE}),
there exist $\ol \tau\in EL(\ol P)$ such that $\ol \tau (\ol q)=(\ol a
p_0,q'')$ and $ht (O_Q(q''))\ol{R[X]}_{\ol a}\geq
\text{rank}(P)-1$. Since $\ol \tau$ can be lifted to $\tau\in
\Aut(P)$, replacing $P$ by $\tau(P)$, we may assume that $ht
(O_Q(q'))\geq \text{rank}(P)-1$ on the Zariski-open set $D(\ol a)$ of
$\Spec(\ol{R[X]})$.

Let $\p_1,\ldots,\p_r$ be minimal prime ideals of $O_Q(q')$ in
$\ol{R[X]}$ not containing $\ol a$. Then $ht(\cap_1^r \p_i)\geq
\text{rank}(P)-1$. By hypothesis $(b)$, we can find $\ol\sigma \in
\Aut(\ol{R[X]})$ with $\ol \sigma(X)=X$ and $\sigma\in \Aut(R[X])$ with
$\sigma(X)=X$ which is a lift of $\ol \sigma$
such that $\ol
\sigma(\cap_1^r \p_i)$ contains a monic polynomial in $\ol {R[X]}=\ol
R[X]$.  Note that $\sigma(f)$ is a monic polynomial. Replacing $R[X]$
by $\sigma(R[X])$, we may assume that $\cap_1^r \p_i$ contains a monic
polynomial in $\ol R[X]$, and $f\in O_P(q)$ is a monic polynomial.

If $\p$ is a minimal prime ideals of $O_Q(q')$ in $\ol{R[X]}$
containing $\ol a$, then $\p$ contains $O_{\ol P}(\ol q)$. Since $f\in
O_P(q)$, $\p$ contains the monic polynomial $\ol f$. Therefore, all
minimal primes of $O_Q(q')$ contains a monic polynomial, hence
$O_Q(q')$ contains a monic polynomial, say $\ol g\in \ol R[X]$. Let
$g\in R[X]$ be a monic polynomial which is a lift of $\ol g$.

{\bf Claim:} For large $N>0$, $p_2=p+X^Ng^Nq \in \Um(P_{1+JR})$.

Choose $\phi \in P^*$ such that $\phi(q)=f$. Then
$\phi(p_2)=\phi(p)+X^Ng^Nf$ is a monic polynomial for large $N$.
Since $p=p_0$ module $JP$, $\ol p=p_0$ and $\ol q=(\ol a \ol
p,q')$. Therefore, 

$\ol p_2=\ol p+X^N\ol g^N (\ol a \ol
p,q')=((1+T^N\ol g^N\ol a)\ol p,X^N\ol g^N q')$.

Since $\ol g \in O_Q(q')\subset O_{\ol P}(\ol p_2)$, we get $O_{\ol
  P}(\ol p) \subset O_{\ol P}(\ol p_2)$. Since $\ol p\in \Um(\ol P)$,
we get $\ol p_2 \in \Um(\ol P)$ and hence $p_2\in \Um(P_{1+JR[X]})$.
Since $O_P(p_2)$ contains a monic polynomial, by \cite[Lemma 1.1, p. 79]{lam}, $p_2\in
\Um(P_{1+JR})$.

Now $p_2=p=p_1$ modulo $XP$, we get $p_2\in \Um(P/XP)$. By (\ref{criterion}),
there exist $p_3\in \Um(P)$ such that $p_3=p_2=p_1$ modulo $XP$. This
completes the proof. $\hfill \square$

\subsection{Proof of Theorem \ref{blr20}} 

Without loss of generality, we may assume that $R$ is reduced.  When
$n=1$, the result follows from well known Quillen \cite{Q} and Suslin
\cite{Su}. When $n=2$, the result follows from Bhatwadekar
\cite[Proposition 3.3]{Bh} where he proves that if $P$ is a projective
$R[T]$-module of rank $2$ such that $P_f$ contains a unimodular
element for some monic polynomial $f\in R[T]$, then $P$ contains a
unimodular element. So now we assume $n\geq 3$.

Write
$A=R[M]$.  Let $J(A,P)=\{s\in A | P_s$ is extended from $A_s\}$ be the
{\it Quillen ideal} of $P$ in $A$.  Let $\wt J=J(A,P)\cap R$ be the ideal
of $R$ and $J=\wt JR[M]$.  We will show that $J$ satisfies the properties of
(\ref{roit1}).

Let $\p\in \Spec(R)$ with $ht(\p)=1$ and $S=R-\p$. Then $S^{-1}P$ is a
projective module over $S^{-1}A[T]=R_{\p}[M][T]$.  Since
$\text{dim}(R_{\p})=1$, by (\ref{ks1}), $S^{-1}P=\wedge^n P_S \oplus
S^{-1}A[T]^{n-1}$. Since determinant of $P$ is extended from $R$,
$\wedge^n P_S=A[T]_S$ and hence $S^{-1}P$ is free. Therefore there
exists $s\in R-\p$ such that $P_s$ is free. Hence $s\in \wt J$
and so $\text{ht}(\wt J)\geq2$.

Since $\text{dim}(R/\wt J)\leq n-2$ and $A[T]/(J)=(R/\wt J)[M][T]$, by
(\ref{ks1}), $P/JP$ contains a unimodular element.

If $I$ is an ideal of $(A/J)[T]=(R/\wt J)[M][T]$ of height $\geq n-1$,
then by (\ref{rem1}), there exists an $R[T]$-automorphism $\sigma\in
\Aut_{R[T]}(A[T])$ such that if $\ol \sigma$ denotes the induced
automorphism of $(A/J)[T]$, then $\ol \sigma (I)$ contains a monic
polynomial in $T$.

By (\ref{dk1}), $EL(P/(T,J)P)$ acts transitively
on $\Um(P/(J,T)P)$.

Therefore, the result now follows from (\ref{roit1}).
\qed
\medskip

\bc\label{blr2}
Let $R$ be a  ring of dimension $n$,
$A=R[X_1,\cdots,X_m]$ a polynomial ring over $R$ and $P$ be a
projective $A[T]$-module of rank $n$.
Assume that $P/TP$ and $P_f$ both contain a unimodular element for some
monic polynomial $f(T)\in A[T]$.  Then $P$ has a unimodular element.
\ec

\proof If $n=1$, the result follows from well known Quillen \cite{Q}
and Suslin \cite{Su} Theorem. When $n=2$, the result follows from
Bhatwadekar \cite[Proposition 3.3]{Bh}. Assume $n\geq 3$. Let $L$ be
the determinant of $P$. If $\wt R$ is the seminormalization of $R$,
then by Swan \cite{Sw}, $L\ot \wt R[X_1,\ldots,X_m]$ is extended from
$\wt R$.  By (\ref{blr20}), $P\ot \wt R[X_1,\ldots,X_m]$ has a
unimodular element.  Since $\wt R[X_1,\ldots,X_n]$ is the
seminormalization of $A$, by Bhatwadekar \cite[Lemma 3.1]{Bh}, $P$ has
a unimodular element.  \qed

\subsection{Obstruction for Projective Modules to have a Unimodular Element}

Let $R$ be a ring of dimension $n\geq 2$ containing $\mathbb{Q}$ and
$P$ be a projective $R[X,Y]$-module of rank $n$ with determinant $L$.
Let $\chi:L\iso \wedge^n(P)$ be an isomorphism. We call $\chi$ an {\it
  orientation} of $P$.  In general, we shall use `hat' when we move to
$R(X)[Y]$ and `bar' when we move modulo the ideal $(X)$. For instance, we have:
\begin{enumerate}
\item $L\ot R(X)[Y]=\hat L$ and $L/XL=\ol L$,
\item $P\ot R(X)[Y]=\hat P$ and $P/XP=\ol P$.
\end{enumerate}
Similarly, $\hat{\chi}$ denotes the induced isomorphism $\hat L\iso
\wedge^n\hat{P}$ and $\ol{\chi}$ denotes the induced isomorphism $\ol
L\iso \wedge^n{\ol P}$.

 We now define the {\it Euler class}  of $(P,\chi)$.  
 
\definition 
First we consider the case $n\geq 2$ and $n\neq 3$.
Let $E(R(X)[Y],\hat L)$ be the $n$th Euler class group of
$R(X)[Y]$ with respect to the line bundle $\hat L$ over $R(X)[Y]$ and
$E(R[Y],\ol L)$ be the $n$th Euler class group of $R[Y]$ with respect
to the line bundle $\ol L$ over $R[Y]$ (see \cite[Section 6]{dz} for
definition). We define the $n$th \emph{ Euler class group} of $R[X,Y]$, 
denoted by $E(R[X,Y],L)$, as the product $E(R(X)[Y],\hat L)\times
E(R[Y],\ol L)$.

 To the pair $(P,\chi)$, 
we associate an element $e(P,\chi)$ of $E(R[X,Y],L)$, called the \emph{Euler class} of $(P,\chi)$, as follows:
$$e(P,\chi)= (e(\hat P, \hat \chi), e(\ol P,\ol \chi))$$
where $e(\hat P, \hat \chi)\in E(R(X)[Y],\hat L)$ is the Euler class of $(\hat P, \hat \chi)$ and $e(\ol P,\ol \chi)\in E(R[Y],\ol L)$ is the 
Euler class of $(\ol P,\ol \chi)$,  defined in \cite[Section 6]{dz}.

Now we treat the case when $n=3$.
Let $\wi E(R(X)[Y],\hat L)$ be the $n$th restricted Euler class group of
$R(X)[Y]$ with respect to the line bundle $\hat L$ over $R(X)[Y]$ and
$\wi E(R[Y],\ol L)$ be the $n$th restricted Euler class group of $R[Y]$ with respect
to the line bundle $\ol L$ over $R[Y]$ (see \cite[Section 7]{dz} for
definition). We define the \emph{Euler class group} of $R[X,Y]$, 
again denoted by $ E(R[X,Y],L)$, as the product $\wi E(R(X)[Y],\hat L)\times
\wi E(R[Y],\ol L)$.

To the pair $(P,\chi)$, 
we associate an element $e(P,\chi)$ of $E(R[X,Y],L)$, called the \emph{Euler class} of $(P,\chi)$, as follows:
$$e(P,\chi)= (e(\hat P, \hat \chi), e(\ol P,\ol \chi))$$
where $e(\hat P, \hat \chi)\in \wt E(R(X)[Y],\hat L)$ is the Euler class of $(\hat P, \hat \chi)$ and $e(\ol P,\ol \chi)\in \wt E(R[Y],\ol L)$ is the 
Euler class of $(\ol P,\ol \chi)$,  defined in \cite[Section 7]{dz}.
\smallskip

\remark
Note that when $n=2$, the definition of the Euler class group $E(R[T], L)$ is slightly different from the case $n\geq 4$.
 See \cite[Remark 7.8]{dz} for details.

\bt\label{obstruct} 
Let $R$ be a ring containing $\mathbb{Q}$ of dimension $n\geq 2$
and $P$ be a projective $R[X,Y]$-module of rank $n$ with determinant
$L$.  Let $\chi:L\iso \wedge^n(P)$ be an isomorphism. Then $e(P,\chi)=
0$ in $E(R[X,Y],L)$ if and only if $P$ has a unimodular element.  \et

\proof First we assume that $P$ has a unimodular element. Therefore,
$\hat P$ and $\ol P$ also have unimodular elements.  If $n\geq 4$, by
\cite[Theorem 6.12]{dz}, we have $e(\hat P, \hat \chi)=0$ in
$E(R(X)[Y],\hat L)$ and $e(\ol P,\ol \chi)=0$ in $E(R[Y],\ol L)$.
The case $n=2$ is taken care by \cite[Remark 7.8]{dz}.  Now if $n=3$,
it follows from \cite[Theorem 7.4]{dz} that $e(\hat P, \hat \chi)=0$
in $ E(R(X)[Y],\hat L)$ and $e(\ol P,\ol \chi)=0$ in $\wi E(R[Y],\ol
L)$.  Consequently, $e(P,\chi)=0$.

Conversely, assume that $e(P,\chi)= 0$. Then $e(\hat P, \hat \chi)=0$ in $E(R(X)[Y],\hat L)$
 and $e(\ol P,\ol \chi)=0$ in $E(R[Y],\ol L)$. 
If $n\neq 3$, by \cite[Theorem 6.12]{dz} and \cite[Remark 7.8]{dz}, $\hat P$ and $\ol P$ have unimodular elements. 
If $n=3$, by \cite[Theorem 7.4]{dz}, $\hat P$ and $\ol P$ have unimodular elements.
Since $\hat P$ has a unimodular element, 
we can find a monic polynomial $f\in R[X]$ such that $P_f$ contains a unimodular element. 
But then by Theorem \ref{blr2},  $P$ has a unimodular element.
\qed
\medskip

\remark Let $R$ be a ring containing $\mathbb{Q}$ of dimension $n\geq 2$ and
$P$ be a projective $R[X_1,\ldots,X_r]$-module ($r\geq 3$) of rank
$n$ with determinant $L$. Let $\chi:L\iso \wedge^r(P)$ be an isomorphism. 
By induction on $r$, we can define the Euler class group of
$R[X_1,\ldots,X_r]$ with respect to the line bundle $L$,
denoted by $E(R[X_1,\ldots,X_r],L)$, as the product of $E(R(X_r)[X_1,\ldots,X_{r-1}],\hat{L})$ and 
$E(R[X_1,\ldots,X_{r-1}],\ol{L}$). 
  
To the pair $(P,\chi)$, we can associate an invariant $e(P,\chi)$ in $E(R[X_1,\ldots,X_r],L)$ as follows:
$$e(P,\chi)= (e(\hat P, \hat \chi), e(\ol P,\ol \chi))$$
where $e(\hat P, \hat \chi)\in  E(R(X_r)[X_1,\ldots,X_{r-1}],\hat L)$ is the Euler class of $(\hat P, \hat \chi)$ and
 $e(\ol P,\ol \chi)\in  E(R[X_1,\ldots,X_{r-1}],\ol L)$ is the Euler class of $(\ol P,\ol \chi)$.
 Finally we have the following result. 

\bt
Let $R$ be a ring containing $\mathbb{Q}$ of dimension $n\geq 2$
and $P$ be a projective $R[X_1,\ldots,X_r]$-module of rank $n$ with determinant
$L$.  Let $\chi:L\iso \wedge^n(P)$ be an isomorphism. Then $e(P,\chi)=
0$ in $E(R[X_1,\ldots,X_r],L)$ if and only if $P$ has a unimodular element.
\et

\section{Analogue of Roy and Mandal}
In this section we will prove (\ref{rm}). We begin with the following
result from \cite[Lemma 2.1]{ry}.

\bl\label{l1}
Let $R$ be a  ring and $P, Q$ be two projective $R$-modules. Suppose that $\phi:Q\longrightarrow P$ is an $R$-linear map.
 For an ideal $I$ of $R$, if $\phi$ is a split monomorphism modulo $I$, then $\phi_{1+I}:Q_{1+I}\longrightarrow P_{1+I}$ is
also a split monomorphism.
\el

\bl Let $(R,\mathcal{M})$ be a local ring and $A$ be a ring such that
$R[X]\hookrightarrow A\hookrightarrow R[X,X^{-1}]$.  Let $P$ and $Q$
be two projective $A$-modules and $\phi:Q\longrightarrow P$ be an
$R$-linear map. If $\phi$ is a split monomorphism modulo $\mathcal{M}$
and if $\phi_f$ is a split monomorphism for some special monic
polynomial $f\in R[X]$, then $\phi$ is also a split monomorphism.  \el 

\proof
By Lemma \ref{l1} $\phi_{1+\mathcal{M}A}$ is a split monomorphism. So,
there is an element $h$ in $1+\mathcal{M}A$ such that $\phi_h$ is a
split monomorphism. Since $f$ is a special monic polynomial,
$R\hookrightarrow A/f$ is an integral extension and hence, $h$ and $f$ are
comaximal.  As $\phi_f$ is also a split monomorphism, it follows that
$\phi$ is a split monomorphism.  \qed 

\bl Let $R$ be a local ring and
$A$ be a ring such that $R[X]\hookrightarrow A\hookrightarrow
R[X,X^{-1}]$.  Let $P$ and $Q$ be two projective $A$-modules and
$\phi,\psi:Q\longrightarrow P$ be $R$-linear maps. Further assume that
$\gamma:P\longrightarrow Q$ is a $A$-linear map such that
$\gamma\psi=f1_Q$ for some special monic polynomial $f\in R[X]$. For large
$m$, there exists a special monic polynomial $g_m\in A$ such that
$X\phi+(1+X^m)\psi$ becomes a split monomorphism after inverting
$g_m$.  \el

\proof
As in \cite{ry, m}, first we assume that $Q$ is free. We have $\gamma(X\phi+(1+X^m)\psi)=X\gamma\phi+(1+X^m)f1_Q$.
Since $Q$ is free, $X\gamma\phi+(1+X^m)f1_Q$ is a matrix. Clearly for large integer $m$, $det(X\gamma\phi+(1+X^m)f1_Q)$ 
 is a special monic polynomial which can be taken for $g_m$.


In the general case, find projective $A$-module $Q'$ such that $Q\op Q'$ is free.
Define maps $\phi', \psi': Q\op Q'\longrightarrow P\op Q'$ and $\gamma': P\op Q'\longrightarrow Q\op Q'$ as
$\phi'=\phi\op 0$, $\psi'=\psi\op f1_{Q'}$ and $\gamma'=\gamma\op 1_{Q'}$. By the previous case, we can find a special monic polynomial $g_m$ 
for some large $m$ such that $(X\phi'+(1+X^m)\psi')_{g_m}$ becomes a split monomorphism. Hence 
$X\phi+(1+X^m)\psi$ becomes a split monomorphism after inverting $g_m$.
\qed

The following result generalizes Mandal's \cite{m}. 
\bt\label{roy}
Let $(R,\mathcal{M})$ be a  local ring and $R[X]\subset A\subset R[X,X^{-1}]$. 
 Let $P$ and $Q$ be two projective $A$-modules with $\text{rank}(Q)< \text{rank}(P)$. If $Q_f$ is a direct summand
of $P_f$ for some special monic polynomial $f\in R[X]$, then $Q$ is also a direct
summand of $P$.
\et

\proof The method of proof is similar to \cite[Theorem 1.1]{ry}, hence
we give an outline of the proof.

Since $Q_f$ is a direct summand of $P_f$,
we can find $A$-linear maps $\psi:Q\longrightarrow P$ and $\gamma :P\longrightarrow Q$ such that $\gamma\psi=f1_Q$
(possibly after replacing $f$ by a power of $f$).

Let 'bar' denote reduction modulo $\mathcal{M}$. Then we have $\bar\gamma\bar\psi=\bar f1_{\bar Q}$. 
As $f$ is special monic, $\bar \psi$ is a monomorphism.

We may assume that $A=R[X,f_1/X^t,\ldots,f_n/X^t]$ with $f_i\in
R[X]$. If $f_i\in \CM R[X]$, then $\ol {R[X,f_i/X^t]}=\ol{R}[X,Y]/(X^tY)$. 
If $f_i\in R[X]-\CM
R[X]$, then $\ol {R[X,f_i/X^t]}$ is either $\ol R[X]$ or $\ol
R[X,X^{-1}]$ depending on whether $\ol f_i/X^t$ is a polynomial in
$\ol R[X]$ or $\ol F_i/X^s$ with $\ol F_i(0)\neq 0$ and $s>0$.

In general, $\ol A$ is one of $\ol R[X]$, $\ol R[X,X^{-1}]$ 
or $\ol R[X,Y_1,\ldots,Y_m]/(X^t(Y_1,\ldots,Y_m))$ for some $m>0$. By \cite[Theorem 3.2]{v}, 
any projective $\ol R[X,Y_1,\ldots,Y_m]/(X^t(Y_1,\ldots,Y_m))$-module is free. 
Therefore, in all cases, projective $\ol A$-modules are free and hence extended
from $\ol R[X]$.  In particular, $\bar P$ and $\bar Q$ are extended
from $\ol R[X]$, which is a PID.
 
Let $\text{rank}(P)=r$ and $\text{rank}(Q)=s$. 
Therefore, using elementary divisors theorem, we can find bases
$\{\bar p_1,\cdots,\bar p_r\}$ and $\{\bar q_1,\cdots,\bar q_s\}$ for
$\bar P$ and $\bar Q$, respectively, such that $\bar \psi(\bar
q_i)=\bar f_i\bar p_i$ for some $f_i\in R[X]$ and $1\leq i\leq s$.

For the rest of the proof, we can follow the proof of
\cite[Theorem 1.1]{ry}.  \qed

\medskip

Now we have the following consequence of (\ref{roy}). 

\bc\label{coro}
Let $R$ be a  local ring and $R[X]\subset A\subset R[X,X^{-1}]$. Let $P$, $Q$ be two
projective $A$-modules such that $P_f$ is isomorphic to $Q_f$ for some special monic polynomial $f\in R[X]$. Then,
\begin{enumerate}
 \item $Q$ is a direct summand of $P\op L$ for any projective $A$-module $L$.
\item $P$ is isomorphic to $Q$ if $P$ or $Q$ has a direct summand of rank one.
\item $P\op L$ is isomorphic to $Q\op L$ for all rank one projective $A$-modules $L$.
\item $P$ and $Q$ have same number of generators.
\end{enumerate}
\ec
\proof
 (1) trivially follows from Theorem \ref{roy} and (3) follows from (2). 

The proof of (4) is same as \cite[Proposition 3.1 (4)]{ry}.

For (2), we can follow the proof of \cite[Theorem 2.2 (ii)]{m} by replacing doubly monic polynomial by special monic polynomial in his arguments.
\qed

\bc
Let $R$ be a  local ring and $R[X]\subset A\subset R[X,X^{-1}]$. Let $P$ be a projective $A$-module such that $P_f$
 is free for some special monic polynomial $f\in R[X]$. Then $P$ is free.
\ec
\proof
Follows from second part of (\ref{coro}).
\qed

\end{document}